\documentclass[a4paper,11pt]{article}
\usepackage{amsfonts}
\usepackage{amssymb}
\usepackage{amsmath}
\usepackage{amsthm}
\usepackage{multirow}
\usepackage{diagbox}
\usepackage{amsopn}
\usepackage{gensymb} 
\usepackage{fancyhdr}
\usepackage{graphicx}
\usepackage{mathrsfs}
\usepackage{latexsym}
\usepackage{graphicx}
\usepackage{subfigure}
\usepackage{color}
 \usepackage{calc}
  \usepackage{tikz}
  
\usepackage{setspace}
\usepackage{textcomp} 
\usepackage{neuralnetwork}

 \usepackage{enumerate}
\usetikzlibrary{decorations.pathreplacing}
 
 \usepackage[noend]{algpseudocode}

\usepackage{algorithmicx,algorithm}
\DeclareMathOperator*{\argmin}{argmin}
 
 \usepackage{appendix}
 
 \usepackage{enumitem}
    \usepackage{libertine}
 
\usetikzlibrary{decorations.pathreplacing}

\numberwithin{equation}{section} \allowdisplaybreaks
\theoremstyle{plain}

\theoremstyle{remark}
\newtheorem{remark}{Remark}
\numberwithin{equation}{section}

\begin{document}


\title{\bf{Bayesian inverse problems using homotopy}}

\author{\hspace{-2.5cm}{\small Xiao-Mei Yang$^1$\thanks{
Corresponding author:  yangxiaomath@swjtu.edu.cn;
Supported by Central Government Funds of Guiding Local Scientific and Technological Development for Sichuan Province No. 2021ZYD0007, NSFC No. 11601067.
} and Zhi-Liang Deng$^{2}$} \\
{\scriptsize $1.$ School of Mathematics,
Southwest Jiaotong University,
Chengdu 610031, China}
\\
{\scriptsize $2.$ School of Mathematical Sciences,  University of Electronic Science and Technology of China,
Chengdu 610054, China}
}
\date{}

\maketitle
\begin{abstract}
\noindent  In solving Bayesian inverse problems, it is often desirable to use a common density parameterization to denote the prior and posterior. 
 Typically we seek a density from the same family as the prior which closely approximates the true posterior. As one of the most important classes of distributions in statistics, the exponential family is considered as the parameterization. The optimal parameter values for representing the approximated posterior are achieved by minimizing the deviation between the parameterized density and a homotopy that deforms the prior density into the posterior density. 
 Rather than trying to solve the original problem, it is exactly converted into a corresponding system
of explicit ordinary first-order differential equations. Solving this system over a finite 'time' interval yields the
desired optimal density parameters.
 This method is proven to be effective by some numerical examples.



%



\noindent \textbf{Key words:}    Inverse scattering; Uncertainty quantification; Bayesian inversion; Homotopy;

\noindent \textbf{MSC 2010}: 35R30, 62F15
\end{abstract}

\section{Introduction}
Inverse problems occur widely in mathematical and engineering fields. The related mathematical theories and algorithms have been developed by many authors \cite{Tarantola}. 
Recently, Bayesian inference establishes a convenient framework in order to analyze the uncertainty of the unknowns \cite{Stuart}. It adopts probability viewpoint to represent, propagate and update epistemic parameter uncertainties.
In computation, a key challenge lies in the transition from the prior distribution to the posterior, which attracts a large number of researchers. 
Some existed algorithms, e.g., the Metropolis-Hastings (MH) sampling algorithm, have been proposed to explore the posterior distribution. Based on MH sampling technique, plenty of samples are collected by reject-accept proposal and are used to characterize the posterior distribution. As we know, the technical difficulty posed by MCMC based methods is that the samples will typically be autocorrelated (or anticorrelated) within a chain. This increases the uncertainty of the estimation of posterior quantities of interest, such as means, variance. 





In this paper, we seek an optimal approximation of the posterior from some common family of distribution. 
We assume that the prior belongs to this distribution class. Therefore, the optimal approximation and the prior are of the same type and share the common parameterization density.
The transition from the prior to the posterior involves the change of the distribution parameters. Homotopy methods are a promising approach to characterize solution spaces by smoothly tracking solutions from one formulation (typically an "easy" problem) to another (typically a ''hard problem"). 
 In deal with filter problems, Hanebeck et al. \cite{Hanebeck} introduced a general framework for performing the measurement update using a homotopy.
This method is discussed further in \cite{Hagmar}. To the authors' knowledge, this method has not been applied to solve inverse problems. We discuss the application of this method in Bayesian inverse problems.
 Different from \cite{Hagmar, Hanebeck}, the approximated posterior density family is taken as the exponential family distribution and the mixed exponential family,  and the usual moment parameters are replaced by a simplition of natural parameters. Within the exponential family and the mixed version, the corresponding derivatives can be computed in a relative easy way. And the dimension number of the parameters decreases dramatically.  With this homotopy, the prior parameters is promoted to the posterior parameters by a so-called homotopy differential equation (HDE). In this HDE, we confront the high-dimension numerical integration. Some conventional means, e.g., sparse grid \cite{Gerstner}, Monte Carlo methods \cite{Dick}, can be used to deal with these integration terms. 
 
 

The remainder of this paper is organized as follows: In Section 2, we give the basic framework of Bayesian inversion using homotopy.  In Section 3, we introduce the exponential family and mixed exponential family. In Section 4, the approximated version of  the homotopy differential equation is derived. In Section 5, some numerical examples are given to verify the effectiveness of the proposed algorithm. 

\section{Bayesian inversion using homotopy}

Inverse problems concern converting observational data into information about systems which are not observed directly. In mathematics, an inverse problem takes the abstract form
\begin{align}\label{s1.1}
\boldsymbol{\tilde{y}}=\mathcal{G}(\boldsymbol{\kappa})+\boldsymbol{\Xi}
\end{align}
in which $\boldsymbol{\Xi}$ is an additive noise,  the unknown $\boldsymbol{\kappa}\in U$ is to be determined, given the data $\boldsymbol{\tilde{y}}\in Y$, where $U$ and $Y$ are Banach spaces. We apply a probabilistic viewpoint, Bayesian approach, to give the solution information of \eqref{s1.1}, in which all
quantities including the unknown $\boldsymbol{\kappa}$, the noise $\Xi$ and the observations $\boldsymbol{\tilde{y}}$ are regarded as random variables. 
In the Bayesian framework, the information about the unknown is updated by blending prior beliefs with observed data. 
Typically, the prior and posterior are coded in the corresponding probability measures, which are linked by the Bayes' formula
\begin{align}
\frac{d\mu^{\boldsymbol{\tilde{y}}}}
{d\mu_0}(\boldsymbol{\kappa})
\propto L(\boldsymbol{\kappa}; \boldsymbol{\tilde{y}}):=\exp\left(-\Phi(\boldsymbol{\kappa}; \boldsymbol{\tilde{y}})\right),
\end{align}
where $L(\boldsymbol{\kappa}; \boldsymbol{\tilde{y}})$ is the likelihood function and $\Phi(\boldsymbol{\kappa}; \boldsymbol{\tilde{y}})$ is the negative log likelihood.
%

For the sake of simplicity, we consider the spaces $U$ and $Y$ are finite dimensional. The posterior density corresponding to $\mu^{\tilde{\boldsymbol{y}}}(\boldsymbol{\kappa})$ is given as
\begin{align}\label{pos1}
\mathfrak{p}^{\boldsymbol{\tilde{y}}}(\boldsymbol{\kappa})
=\frac{\exp\left(-\Phi(\boldsymbol{\kappa}; \boldsymbol{\tilde{y}})\right)\mathfrak{q}(\boldsymbol{\kappa})}{Z}.
\end{align}
The denominator $Z$ is the normalization constant, which is usually neglected in sampling algorithms. The main target to be explored in the posterior density is the numerator 
\begin{align}
\mathfrak{p}(\boldsymbol{\kappa})=\exp\left(-\Phi(\boldsymbol{\kappa}; \boldsymbol{\tilde{y}})\right)\mathfrak{q}(\boldsymbol{\kappa}),
\end{align}
where the normalization constant $Z$ in \eqref{pos1} is dropped. We try to find an optimal approximation of $\mathfrak{p}(\boldsymbol{\kappa})$ from some commonly used distribution family.  This is implemented in the frame of homotopy Bayesian approach.


The key idea of the homotopy Bayesian approach is to perform progressive processing.
 In this method, instead of directly approximating the true density $\mathfrak{p}(\boldsymbol{\kappa})$, it starts with a tractable density and continuously approaches the true density via intermediate densities. 
 Choose the homotopy as follows
\begin{align}
\mathfrak{p}(\boldsymbol{\kappa}, t)=L^t(\boldsymbol{\kappa}; \boldsymbol{\tilde{y}})\mathfrak{q}(\boldsymbol{\kappa})=\exp\left(-t\Phi(\boldsymbol{\kappa}; \boldsymbol{\tilde{y}})\right)\mathfrak{q}(\boldsymbol{\kappa}), \,\, t\in [0, 1].
\end{align}
For each $t\in [0, 1]$, an optimal approximation of $\mathfrak{p}(\boldsymbol{\kappa}, t)$ is to be sought from some distribution family. Denote the probability density of the distribution family by a parameterized density family $\mathfrak{g}(\boldsymbol{\kappa}; \eta)$. 
By miminzing a deviation function 
 $G(\eta, t)$ between $\mathfrak{p}(\boldsymbol{\kappa}, t)$ and  $\mathfrak{g}(\boldsymbol{\kappa}; \eta)$,  we obtain the optimal approximation of $\mathfrak{p}(\boldsymbol{\kappa}, t)$. 

To measure the difference between two probability distributions over the same variable $\boldsymbol{\kappa}$, a measure, called the Kullback-Leibler divergence, or simply, the KL divergence, has been popularly used in the statistical learning, data mining literature. The concept was originated in probability theory and information theory. The KL divergence, which is closely related to relative entropy, information divergence, and information for discrimination, is a non-symmetric measure of the difference between two probability distributions $\mathfrak{X}$ and $\mathfrak{Y}$. If $\mathfrak{X}$ and $\mathfrak{Y}$ are discrete probability distributions, i.e., $\mathfrak{X}=(\mathfrak{X}(\boldsymbol{\kappa}_1), \mathfrak{X}(\boldsymbol{\kappa}_2), \cdots, \mathfrak{X}(\boldsymbol{\kappa}_m))$ and $\mathfrak{Y}=(\mathfrak{Y}(\boldsymbol{\kappa}_1), \mathfrak{Y}(\boldsymbol{\kappa}_2), \cdots, \mathfrak{Y}(\boldsymbol{\kappa}_m))$, the KL divergence is defined to be
\begin{align}
D_{\rm KL}(\mathfrak{X}||\mathfrak{Y})=\sum_{i=1}^m\mathfrak{X}(\boldsymbol{\kappa}_i)\log\frac{\mathfrak{X}(\boldsymbol{\kappa}_i)}{\mathfrak{Y}(\boldsymbol{\kappa}_i)}.
\end{align}
For continuous random variables $\mathfrak{X}$ and $\mathfrak{Y}$, we assume the corresponding probability densities are $\mathfrak{x}(\boldsymbol{\kappa})$ and $\mathfrak{y}(\boldsymbol{\kappa})$, respectively. 
In this case, the  Kullback-Leibler divergence between $\mathfrak{X}$ and $\mathfrak{Y}$ is given
\begin{align}
D_{\rm KL}(\mathfrak{X}||\mathfrak{Y})=\int \mathfrak{x}(\boldsymbol{\kappa})\log\frac{\mathfrak{x}(\boldsymbol{\kappa})}{\mathfrak{y}(\boldsymbol{\kappa})}d\boldsymbol{\kappa}.
\end{align}
In probability and statistics, the Hellinger distance is also usually used to quantify the similarity between two probability distributions. It is a type of f-divergence. 
The squared Hellinger distances for discrete and continuous random variables are defined by
\begin{align}
\begin{aligned}
&H^2(\mathfrak{X}, \mathfrak{Y})=\frac{1}{2}\sum_{i=1}^m \left(\sqrt{\mathfrak{X}(\boldsymbol{\kappa}_i)}-\sqrt{\mathfrak{Y}(\boldsymbol{\kappa}_i)}\right)^2,\\
&H^{2}(\mathfrak{X}, \mathfrak{Y})=\frac{1}{2} \int\left(\sqrt{\mathfrak{x}(\boldsymbol{\kappa})}-\sqrt{\mathfrak{y}(\boldsymbol{\kappa})}\right)^{2} d \boldsymbol{\kappa},
\end{aligned}
\end{align}
respectively. 

With the preceding two measure deviations acting as $G$, we seek
\begin{align}
\eta(t)=\argmin_{\eta} G(\eta, t).
\end{align}
The minimization necessary condition yields
\begin{align}
G_\eta(\eta(t), t)=0 \,\, \text{for} \,\,  t\in [0, 1].
\end{align}
By taking the total derivative w.r.t. $t$ of the preceding equation, we obtain
\begin{align}
G_{\eta\eta}(\eta(t), t)\eta'(t)+G_{\eta t}(\eta(t), t)=0.
\end{align}
For $t=0$, the approximated density $\mathfrak{g}(\boldsymbol{\kappa}; \eta(0))$ is an optimal approximation $\mathfrak{q}(\boldsymbol{\kappa})$. 
If the approximated distributions $\mathfrak{g}(\boldsymbol{\kappa}; \eta)$ are chosen from the same family as the prior distribution $\mathfrak{q}(\boldsymbol{\kappa})$,  the initial condition is set according to $\mathfrak{g}(\boldsymbol{\kappa}; \eta(0))=\mathfrak{q}(\boldsymbol{\kappa})$. Let $\mathfrak{q}(\boldsymbol{\kappa})$ be identified as $\mathfrak{g}(\boldsymbol{\kappa}; \eta_0)$.  By this, the initial condition is set to $\eta(0)=\eta_0$.
Solving for $\eta'(t)$, we arrive at the initial value problem of the homotopy differential equation (HDE)
\begin{align}\label{hde1}
\begin{aligned}
&\eta'(t)=-G_{\eta\eta}^{-1}(\eta(t), t)G_{\eta t}(\eta(t), t)),\\
&\eta(0)=\eta_0.
\end{aligned}
\end{align}


%



%
%
%

\section{Exponential family and mixed exponential family}
%

As the most widely used distribution family, the exponential family (EF) and the mixed exponential family (MEF) are served as the approximation of the posterior distribution.
The EF and MEF are a practically convenient and widely used unified families of distributions parametrized by a finite dimensional parameter vector.
The reason of its special significance
  is that a number of important and useful calculations in statistics
can be done all at one stroke within the framework of the EF and MEF. This generality
contributes to both convenience and larger scale understanding. And besides,  it has recently obtained additional importance due to its use and appeal to the machine learning community.

For a numeric random variable $\boldsymbol{\mathfrak{K}}$, the parametric EF probability density can be written as
\begin{align}
\mathfrak{q}(\boldsymbol{\kappa}; \theta)=h(\boldsymbol{\kappa}) \exp \left\{\langle T(\boldsymbol{\kappa}), \theta\rangle-A(\theta)\right\}
\end{align}
where  $\theta$ is called the natural (canonical) parameter, $T(\boldsymbol{\kappa})$ is the sufficient statistic and $A(\theta)$  is the $\log$ normalizer
given by
\begin{align*}
A(\theta)=\log\int h(\boldsymbol{\kappa})\exp(\langle T(\boldsymbol{\kappa}), \theta\rangle)d\boldsymbol{\kappa}=\log Q(\theta).
\end{align*}
 It is easy to know that
\begin{align}
\frac{d\log \mathfrak{q}(\boldsymbol{\kappa}; \theta)}{d\theta}=T(\boldsymbol{\kappa})-A'(\theta).
\end{align}
We can compute $A'(\theta)$ by 
\begin{align*}
\begin{aligned}
A'(\theta)&=\frac{1}{Q(\theta)}\frac{dQ(\theta)}{d\theta}=\frac{\int h(\boldsymbol{\kappa})\exp(\langle T(\boldsymbol{\kappa}), \theta\rangle)T(\boldsymbol{\kappa})d\boldsymbol{\kappa}}{\int h(\boldsymbol{\kappa})\exp(\langle T(\boldsymbol{\kappa}), \theta\rangle)d\boldsymbol{\kappa}}\\
&=\frac{\int h(\boldsymbol{\kappa})\exp(\langle T(\boldsymbol{\kappa}), \theta\rangle-A(\theta))T(\boldsymbol{\kappa})d\boldsymbol{\kappa}}{\int h(\boldsymbol{\kappa})\exp(\langle T(\boldsymbol{\kappa}), \theta\rangle-A(\theta))d\boldsymbol{\kappa}}\\
&=\mathbb{E}[T(\boldsymbol{\kappa})].
\end{aligned}
\end{align*}

We only list the case of the multivariate Gaussian distribution, which is used in this paper. For this case, we further reduce the parameter dimension number. 
For a Gaussian random variable $\boldsymbol{\mathfrak{K}}\in \mathbb{R}^{d}$, if $\boldsymbol{\mathfrak{K}}\sim N(\varkappa, \varSigma)$, then $\mathbb{E}[\boldsymbol{\mathfrak{K}}]=\varkappa$ 
and $\text{cov}[\boldsymbol{\mathfrak{K}}]=\varSigma$. $\varkappa$ and $\varSigma$ are called the moment parameters of the distribution. 
The probability density is given
\begin{align*}
&\mathfrak{q}(\boldsymbol{\kappa} \mid \varkappa, \varSigma) =\frac{1}{(2 \pi)^{d / 2} \mid \varSigma\mid^{1 / 2}} \exp \left\{-\frac{1}{2}(\boldsymbol{\kappa}-\varkappa)^{\top} \varSigma^{-1}(\boldsymbol{\kappa}-\varkappa)\right\}
\\&=\frac{1}{(2 \pi)^{d/ 2}} \exp \left\{-\frac{1}{2} \operatorname{tr}\left(\varSigma^{-1} \boldsymbol{\kappa} \boldsymbol{\kappa}^{\top}\right)+\varkappa^{\top} \varSigma^{-1} \boldsymbol{\kappa}-\frac{1}{2} \varkappa^{\top} \varSigma^{-1} \varkappa-\frac{1}{2}\log |\varSigma|\right\}.
\end{align*}
The corresponding function $T(\boldsymbol{\kappa})=\left[\begin{array}{c}\boldsymbol{\kappa} \\ \operatorname{vec}(\boldsymbol{\kappa}\boldsymbol{\kappa}^\top)\end{array}\right]$ and natural parameter is
\begin{align*}
\theta=\left[\begin{array}{c}\varSigma^{-1}\varkappa \\ -\frac{1}{2}\operatorname{vec}(\varSigma^{-1})\end{array}\right].
\end{align*} 
Here $A(\theta)=\frac{1}{2} \varkappa^{\top} \varSigma^{-1} \varkappa+\frac{1}{2}\log |\varSigma|$
 and $h(\boldsymbol{\kappa})=(2\pi)^{-\frac{d}{2}}$. Define the precision matrix by $\mathcal{P}=\varSigma^{-1}$. By the symmetric positive definiteness of $\mathcal{P}$,  we have the  Cholesky factorization $\mathcal{P}=\mathcal{R}^\top \mathcal{R}$ with $\mathcal{R}$ being a lower triangular matrix. Introducing the new parameter
 \begin{align}\label{par1}
 \theta=\left[\begin{array}{c}\varkappa \\ \operatorname{vech}(\mathcal{R})\end{array}\right]=\left[\begin{array}{c}\theta_1 \\ \theta_2\end{array}\right],
 \end{align}
we have
\begin{align}
&\frac{\partial \log\mathfrak{q}(\boldsymbol{\kappa} \mid \varkappa, \varSigma) }{\partial \theta_1}=-\mathcal{R}^\top\mathcal{R}(\boldsymbol{\kappa}-\varkappa),\\
&\frac{\partial \log\mathfrak{q}(\boldsymbol{\kappa} \mid \varkappa, \varSigma) }{\partial \theta_2}=\operatorname{vech}\left([\operatorname{diag}(\mathcal{R})]^{-1}-\mathcal{R}(\boldsymbol{\kappa}-\varkappa)(\boldsymbol{\kappa}-\varkappa)^\top\right).
\end{align}
\begin{remark}
It is obvious that the dimension number of the new parameter system decrease dramatically compared with the moment parameter or the natural parameter system. 
\end{remark}

The MEF is denoted by
\begin{align}
\mathfrak{g}(\boldsymbol{\kappa}; \eta)=\sum_{i=1}^{M} w_i \mathfrak{q}(\boldsymbol{\kappa}; \theta^i)=\sum_{i=1}^M w_i h(\boldsymbol{\kappa})\exp\left(\langle \Phi(\boldsymbol{\kappa}), \theta^i \rangle-A_i(\theta^i)\right),
\end{align}
where $\sum_{i=1}^M w_i=1$ and $w_i\geq 0$ for $i=1, 2, \cdots, M$. To remove the weight constraint, we let $w_i=\frac{\pi/2+\arctan\lambda_i}{\sum_{j=1}^M  \pi/2+\arctan\lambda_j}$ and $\lambda_M=0$.
Thus, the MEF has the following form
\begin{align*}
\mathfrak{g}(\boldsymbol{\kappa}; \eta)=\sum_{i=1}^M \frac{\pi/2+\arctan\lambda_i}{\sum_{j=1}^M  \pi/2+\arctan\lambda_j}h(\boldsymbol{\kappa})\exp\left(\langle \Phi(\boldsymbol{\kappa}), \theta^i \rangle-A_i(\theta^i)\right).
\end{align*}
Then it follows by simple calculations that
\begin{align}
&\frac{\partial \log\mathfrak{g}(\boldsymbol{\kappa}; \eta) }{\partial \theta^i}=\frac{w_i \mathfrak{q}(\boldsymbol{\kappa}; \theta^i)}{\mathfrak{g}(\boldsymbol{\kappa}; \eta)}\frac{\partial \log\mathfrak{q}(\boldsymbol{\kappa}; \theta^i) }{\partial \theta^i},\\
&\frac{\partial \log\mathfrak{g}(\boldsymbol{\kappa}; \eta) }{\partial \lambda_i}=\frac{\frac{  \mathfrak{q}(\boldsymbol{\kappa}; \theta^i)}{\mathfrak{g}(\boldsymbol{\kappa}; \eta)}-1}{(1+\lambda_i^2)\sum_{j=1}^M \left( \pi/2+\arctan\lambda_j\right)}, \,\, i=1, 2, \cdots, M-1. 
\end{align}


\begin{remark}
The parameter dimension number of the mixed Gaussian exponential family is $M-1+M(d+d(d+1)/2)$.
\end{remark}

\section{Approximation of HDE}
In HDE \eqref{hde1}, the posterior density is involved in the integration term. It is obvious that the corresponding numerical evaluations become complicated. 
Therefore, it is necessary to give an approximation version of \eqref{hde1} in real numerical simulation process. 


 For a partition $0=t_0<t_1<\cdots<t_n=1$, we have from Taylor expansion of $\eta(t)$ 
\begin{align}
\eta(t_{i+1})=\eta(t_i)+\eta'(t_i)\Delta t_i+o(\Delta t_i), \,\, i=0, \cdots, n-1, 
\end{align}
where $\Delta t_i=t_{i+1}-t_i$. Truncating the Taylor expansion after the linear term, we obtain 
\begin{align}\label{tay1}
\eta_{i+1}=\eta_{i}+\eta'_i \Delta t_i,
\end{align}
where $\eta_i$ and $\eta'_i$ are the approximations of $\eta(t_i)$ and $\eta'(t_i)$ respectively. As an approximation of $\eta'(t_i)$, we can take 
\begin{align}\label{hdel1}
\eta'_i=-\tilde{G}_{\eta\eta}^{-1}(\eta_i, t_i)\tilde{G}_{\eta t}(\eta_i, t_i),
\end{align}
where $\tilde{G}_{\eta\eta}(\eta_i, t_i)$ and $\tilde{G}_{\eta t}(\eta_i, t_i)$ are some approximations of $G_{\eta\eta}(\eta(t_i), t_i)$ and $G_{\eta t}(\eta(t_i), t_i)$ respectively. 

Let us check the two deviations of measures, i.e., the Kullback-Leibler divergence and squared Hellinger metric.
For Kullback-Leibler divergence, we have 
\begin{align}
G(\eta, t)=\int \mathfrak{p}(\boldsymbol{\kappa}, t)\log\frac{\mathfrak{p}(\boldsymbol{\kappa}, t)}{\mathfrak{g}(\boldsymbol{\kappa}; \eta)}d\boldsymbol{\kappa}.
\end{align}
It follows that 
\begin{align}
&G_\eta=-\int \mathfrak{p}(\boldsymbol{\kappa}, t)\frac{\partial \log \mathfrak{g}(\boldsymbol{\kappa}; \eta)}{\partial \eta}d\boldsymbol{\kappa},\\
&
\begin{aligned}
G_{\eta\eta}&=-\int \mathfrak{p}(\boldsymbol{\kappa}, t)\frac{\partial^2 \log \mathfrak{g}(\boldsymbol{\kappa}; \eta)}{\partial \eta^2}d\boldsymbol{\kappa}\\
&=-\int \mathfrak{p}(\boldsymbol{\kappa}, t) \frac{\mathfrak{g}(\boldsymbol{\kappa}; \eta)\frac{\partial^2 \mathfrak{g}(\boldsymbol{\kappa}; \eta)}{\partial\eta^2}-\frac{\partial \mathfrak{g}}{\partial\eta}\frac{\partial \mathfrak{g}}{\partial\eta}^\top}{\mathfrak{g}(\boldsymbol{\kappa}; \eta)^2}d\boldsymbol{\kappa},
\end{aligned}
\\
&
\begin{aligned}
G_{\eta t}&=\int \Phi(\boldsymbol{\kappa}; \boldsymbol{\tilde{y}}) \mathfrak{p}(\boldsymbol{\kappa}, t)\frac{\partial \log \mathfrak{g}(\boldsymbol{\kappa}; \eta)}{\partial \eta}d\boldsymbol{\kappa}\\
&=\int \Phi(\boldsymbol{\kappa}; \boldsymbol{\tilde{y}}) \mathfrak{p}(\boldsymbol{\kappa}, t)\frac{\frac{\partial \mathfrak{g}(\boldsymbol{\kappa}; \eta)}{\partial \eta}}{\mathfrak{g}(\boldsymbol{\kappa}; \eta)}d\boldsymbol{\kappa}.
\end{aligned}
\end{align}
When we implement the iteration \eqref{tay1} with \eqref{hdel1}, the density $\mathfrak{p}(\boldsymbol{\kappa}, t)$ at $t_i$ has been approximated by $\mathfrak{g}(\boldsymbol{\kappa}; \eta_i)$. 
Therefore, we can take the approximations of $G_{\eta\eta}(\eta(t_i), t_i)$ and $G_{\eta t}(\eta(t_i), t_i)$ as 
\begin{align}\label{g11}
\begin{aligned}
\tilde{G}_{\eta\eta}(\eta_i, t_i)&=-\int  \frac{\mathfrak{g}(\boldsymbol{\kappa}; \eta)\frac{\partial^2 \mathfrak{g}(\boldsymbol{\kappa}; \eta)}{\partial\eta^2}-\frac{\partial \mathfrak{g}}{\partial\eta}\frac{\partial \mathfrak{g}}{\partial\eta}^\top}{\mathfrak{g}(\boldsymbol{\kappa}; \eta)}\bigg |_{\eta_i}d\boldsymbol{\kappa}\\
&=\int \frac{\partial \log\mathfrak{g}}{\partial\eta} \frac{\partial \log\mathfrak{g}}{\partial\eta} ^\top \mathfrak{g} \bigg |_{\eta_i}d\boldsymbol{\kappa}-\int \frac{\partial^2\mathfrak{g}}{\partial\eta^2} \bigg |_{\eta_i}d\boldsymbol{\kappa}\\
&=\int \frac{\partial \log\mathfrak{g}}{\partial\eta} \frac{\partial \log\mathfrak{g}}{\partial\eta} ^\top \mathfrak{g} \bigg |_{\eta_i}d\boldsymbol{\kappa}:=I(\eta_i)
\end{aligned}
\end{align}
and
\begin{align}
\begin{aligned}
&\tilde{G}_{\eta t}(\eta_i, t_i)=\int \Phi(\boldsymbol{\kappa}; \boldsymbol{\tilde{y}}) \frac{\partial \mathfrak{g}(\boldsymbol{\kappa}; \eta)}{\partial \eta}\bigg |_{\eta_i}d\boldsymbol{\kappa}.
\end{aligned}
\end{align}
In \eqref{g11}, the term $I(\eta)$ is the Fisher information matrix when $\mathfrak{g}(\boldsymbol{\kappa}; \eta)$ is a probability density.

For squared Hellinger metric, we obtain the corresponding expressions
\begin{align}
&\tilde{G}_{\eta\eta}(\eta_i, t_i)=\frac{1}{4}I(\eta_i),\\
&\tilde{G}_{\eta t}(\eta_i, t_i)=\frac{1}{4}\int \Phi(\boldsymbol{\kappa}; \boldsymbol{\tilde{y}}) \frac{\partial \mathfrak{g}(\boldsymbol{\kappa}; \eta)}{\partial \eta}\bigg |_{\eta_i}d\boldsymbol{\kappa}.
\end{align}

Therefore, for both cases, the approximated homotopy difference equation is
\begin{align}\label{hde4.1}
\eta_{i+1}=\eta_i-\eta'_i\Delta t_i,
\end{align}
where $\eta'_i=I(\eta_i)^{-1}\int \Phi(\boldsymbol{\kappa}; \boldsymbol{\tilde{y}}) \frac{\partial \mathfrak{g}(\boldsymbol{\kappa}; \eta)}{\partial \eta}\mid_{\eta_i}d\boldsymbol{\kappa}$. The corresponding approximated homotopy differential equation can be written as
\begin{align}
\frac{d\tilde{\eta}}{dt}=-I^{-1}(\tilde{\eta})\int \Phi(\boldsymbol{\kappa}; \boldsymbol{\tilde{y}})\frac{\partial\mathfrak{g}}{\partial\tilde{\eta}}d\boldsymbol{\kappa}.
\end{align}


\section{Numerical examples}

The homotopy Bayes algorithm is easy to implement, which we list in Algorithm \ref{alg1} 
\begin{algorithm}
\caption{Homotopy Bayes algorithm}
\begin{algorithmic}[1]

    \State Initial condition: Take the initial parameter from the prior distribution. 
    \State While $t=t_i$,  compute the integration in \eqref{hde4.1} using the Monte Carlo algorithm. And set $\eta_{i+1}=\eta_i-\eta'_i\Delta t_i$ and set $t_{i+1}=t_i+\Delta t$. 
    \State Stop when $t=1$.
   

\end{algorithmic}
\label{alg1}
\end{algorithm}

\subsection{Inverse heat conduction problem}
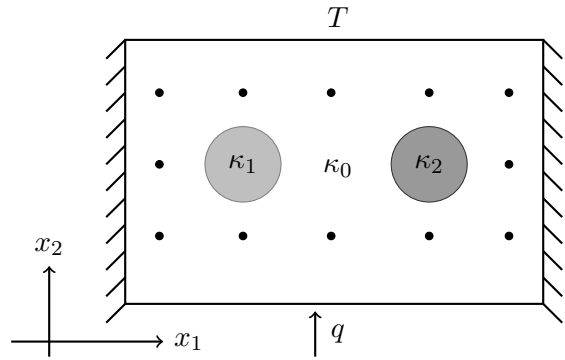
\begin{figure}
\centering
\begin{tikzpicture}
\draw[->, thick] (-3,-2)--(-1,-2) node[right]{$x_1$};
\draw[->, thick] (-2.5,-2.2)--(-2.5,-1) node[above]{$x_2$};
\draw[thick] (-1.5,-1.5) rectangle (4,2);
\draw [thick] (-1.5,-1.5) -- (-1.75,-1.75); 
\draw [thick] (-1.5,-1.15) -- (-1.75,-1.4); 
\draw [thick] (-1.5,-0.8) -- (-1.75,-1.05); 
\draw [thick] (-1.5,-0.45) -- (-1.75,-0.7); 
\draw [thick] (-1.5,-0.1) -- (-1.75,-0.35); 
\draw [thick] (-1.5,0.25) -- (-1.75,0); 
\draw [thick] (-1.5,0.6) -- (-1.75,0.35); 
\draw [thick] (-1.5,0.95) -- (-1.75,0.7); 
\draw [thick] (-1.5,1.3) -- (-1.75,1.05); 
\draw [thick] (-1.5,1.65) -- (-1.75,1.4);
\draw [thick] (-1.5,2) -- (-1.75,1.75);

\draw [thick] (4,2) -- (4.25,1.75); 
\draw [thick] (4,1.65) -- (4.25,1.4); 
\draw [thick] (4,1.3) -- (4.25,1.05); 
\draw [thick] (4,0.95) -- (4.25,0.7);
\draw [thick] (4,0.6) -- (4.25,0.35);  
\draw [thick] (4,0.25) -- (4.25,0); 
\draw [thick] (4,-0.1) -- (4.25,-0.35); 
\draw [thick] (4,-0.45) -- (4.25,-0.7); 
\draw [thick] (4,-0.8) -- (4.25,-1.05); 
\draw [thick] (4,-1.15) -- (4.25,-1.4); 
\draw [thick] (4,-1.5) -- (4.25,-1.75); 

\draw [black,fill=gray, ,opacity=0.5](0.05,0.35) circle (0.5); 

\draw [black,fill=gray, ,opacity=0.8](2.5,0.35) circle (0.5);

\node at (0.05,0.35) {$\kappa_1$};
\node at (2.5,0.35) {$\kappa_2$};
\node at (1.3,0.3) {$\kappa_0$};

\node at (1.3,2.3) {$T$};

\draw[->,  thick] (1,-2.2)--(1,-1.6);

\node at (1.3,-1.9) {$q$};


\foreach \i in {-1.05, 0.05, 1.21, 2.5, 3.55} {\draw (\i, 1.3) [fill=black] circle (0.05);}

\foreach \i in {-1.05, 0.05, 1.21, 2.5, 3.55} {\draw (\i, -0.6) [fill=black] circle (0.05);}

\foreach \i in {-1.05, 3.55} {\draw (\i, 0.35) [fill=black] circle (0.05);}

\end{tikzpicture}

\caption{ 2D IHCP: heat conduction setup.}
\label{fig1}
\end{figure}

We consider the reconstruction of $\boldsymbol{\kappa}$ using the measurements of solution $u$ governed by the following system
\begin{align}\label{eqn1}
\begin{aligned}
&\nabla\cdot(\boldsymbol{\kappa}\nabla u)=0 \,\,\text{in}\,\, [0, 1]\times[0, 0.6],\\
&u(0, x_2)=u(1, x_2)=0, \,\, u(x_1, 0.6)=T=200, \,\, -\kappa_0\frac{\partial u}{\partial x_2}=2000.
\end{aligned}
\end{align}
This example is originally presented in Nagel and Sudret \cite{Nagel} (also \cite{Wagner}). We again solve the same problem with homotopy Bayesian approach and investigate the performance of the algorithm. 
As displayed in Fig. \ref{fig1}, the background thermal conductivity  is denoted as $\kappa_0$ that is known, while the conductivities of the material inclusions are termed as $\kappa_1$ and $\kappa_2$, respectively. 
 
 We consider the inverse heat conduction problem (IHCP) that is posed when the thermal conductivities ${\boldsymbol{\kappa}}=(\kappa_1, \kappa_2)^\top$ are unknown and their inference is intended. 
 With this in mind, a number of $N$ measurements ${\boldsymbol{\tilde{y}}}=(u(\boldsymbol{x}_1), \cdots, u(\boldsymbol{x}_N))^\top$ of the temperature field at the measurement locations $(\boldsymbol{x}_1, \boldsymbol{x}_2, \cdots, \boldsymbol{x}_N)^\top$  is available. The forward model $\mathcal{G}: \boldsymbol{\kappa}\rightarrow \tilde{\boldsymbol{y}}$ is established by the finite element discretization for \eqref{eqn1}. With the discretization model, the measured temperatures $\boldsymbol{\tilde{y}}$ are  generated by adding ansatz noise to the numerical solutions as follows
 \begin{align}
 \boldsymbol{\tilde{y}}=\boldsymbol{y}+\Xi=\mathcal{G}(\boldsymbol{\kappa})+\Xi,
 \end{align}
 where $\Xi\sim N(0, \delta^2 I)$. The prior is set to a multivariate lognormal distribution $\mathfrak{q}(\boldsymbol{\kappa})=\prod_{i=1}^2 \mathfrak{q}(\kappa_i)$ with independent marginals $\mathfrak{q}(\kappa_i)=LN(\kappa_i|\varkappa_0, \sigma_0)$ with $\varkappa_0=30$ and $\sigma_0=6$. 
 Parameters $\varkappa_0$ and $\sigma_0$ describe the mean and standard deviation of the lognormal prior. They are related to the parameters of the associated normal distribution $N(\log(\kappa_i)|\lambda_0, \zeta_0^2)$
  via $\varkappa_0=\exp(\lambda_0+\zeta_0^2/2)$ and $\sigma_0^2=(\exp(\zeta_0^2)-1)\exp(2\lambda_0+\zeta_0^2)$.
 The unknown parameters are represented as $\kappa_i=\exp(\lambda_0+\zeta_0\xi_i)$  in terms of the standardized variables $\xi_i\in\mathbb{R}$ with Gaussian weight functions $N(\xi_i|0, 1)$. In this computation process, we take a uniform partition for homotopy parameter interval 
 $[0, 1]$. In this 2D IHCP, we take the  $\Delta t_i=0.01$. This integrations in \eqref{hde4.1} are computed by the Monte Carlo numerical integration. For noise $\delta=0.25$, we display the numerical results in Fig. \ref{fig3}. The left one in Fig. \ref{fig3} is the function $\kappa=\exp(\lambda_0+\zeta_0\xi)$, where $\xi$ is estimated by the proposed algorithm. The right one in Fig. \ref{fig3} shows the true posterior density.
 
\begin{figure}
\centering
\begin{tikzpicture}
\draw[thick] (0,0) rectangle (8,4);
\foreach \i in {1.8, 4, 6.191}
\foreach \j in {1.2, 3}
{
\fill (\i, \j)[black,fill=gray, ,opacity=(\i+\j)/10] circle (0.5cm);
}

\node at (1.8, 1.2) {$\kappa_{4}$};
\node at (1.8, 3) {$\kappa_{1}$};
\node at (4, 1.2) {$\kappa_{5}$};
\node at (4, 3) {$\kappa_{2}$};
\node at (6.191, 1.2) {$\kappa_{6}$};
\node at (6.191, 3) {$\kappa_{3}$};
\node at (4, 2) {$\kappa_0$};

\foreach \x in {0.7, 1.8, 2.9, 4, 5.1, 6.2, 7.4}
\fill[color=black] (\x, .3) circle (.06cm);

\foreach \x in {0.7, 1.8, 2.9, 4, 5.1, 6.2, 7.4}
\fill[color=black] (\x, 3.7) circle (.06cm);

\foreach \y in {1.1, 2, 2.9}
\fill[color=black] (0.7,\y) circle (.06cm);

\foreach \y in {1.1, 2, 2.9}
\fill[color=black] (7.4,\y) circle (.06cm);

\foreach \y in {0, 0.2, 0.4,0.6,0.8,1,1.2,1.4,1.6,1.8,2,2.2,2.4,2.6,2.8,3,3.2,3.4,3.6,3.8,4}
\draw [thick] (0,\y) -- (-0.2,\y-0.1);

\foreach \y in {0, 0.2, 0.4,0.6,0.8,1,1.2,1.4,1.6,1.8,2,2.2,2.4,2.6,2.8,3,3.2,3.4,3.6,3.8,4}
\draw [thick] (8,\y) -- (8.2,\y-0.1);

\node at (4,4.4) {$T$};

\draw[->,  thick] (4,-1)--(4,-0.2);

\node at (4.3,-0.6) {$q$};

\end{tikzpicture}
\caption{6D IHCP: heat conduction setup.}
\label{fig2}
\end{figure}
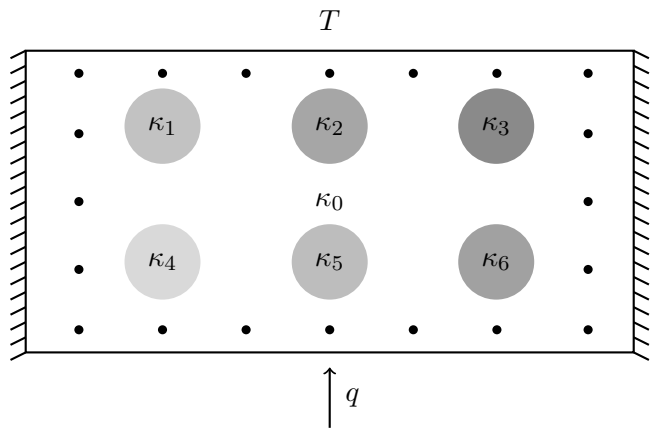

\begin{figure}

\centering

\includegraphics[scale=0.25]{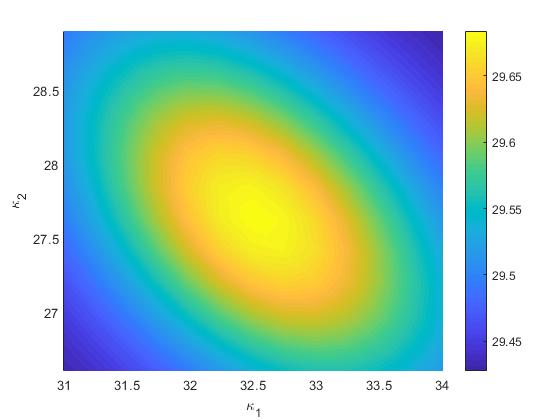}
\includegraphics[scale=0.25]{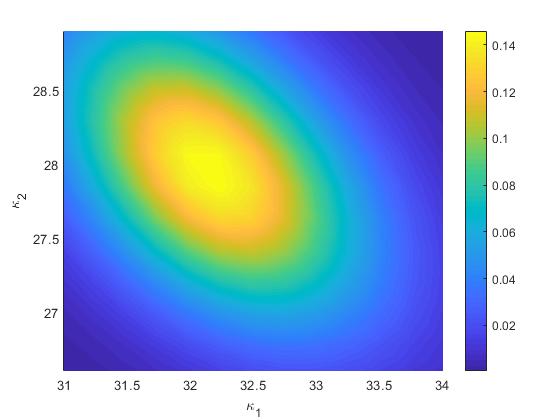}\\
\includegraphics[scale=0.25]{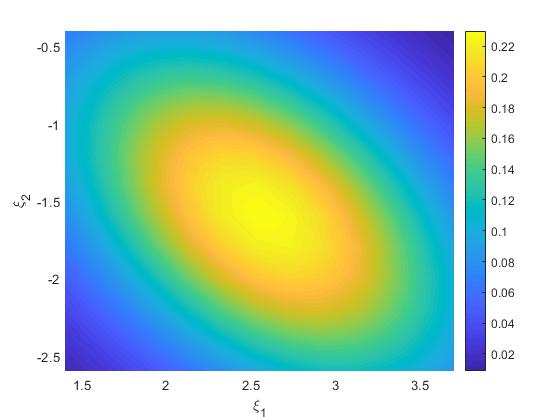}
\includegraphics[scale=0.25]{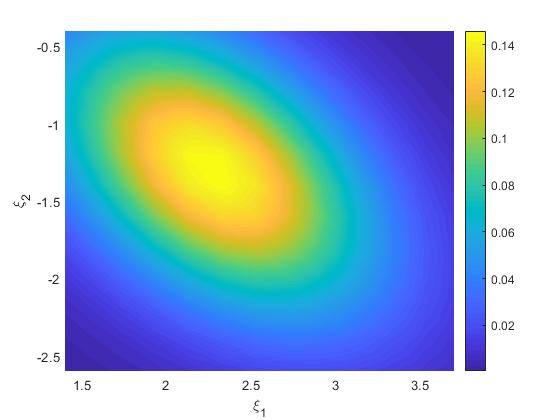}
\caption{The numerical results for the 2D IHCP: The reconstruction (left) for $\boldsymbol{\kappa}=\exp(\lambda_0+\zeta_0\boldsymbol{\xi})$ and the true posterior density (right).}
\label{fig3}

\end{figure}

\begin{figure}

\centering

\includegraphics[scale=0.25]{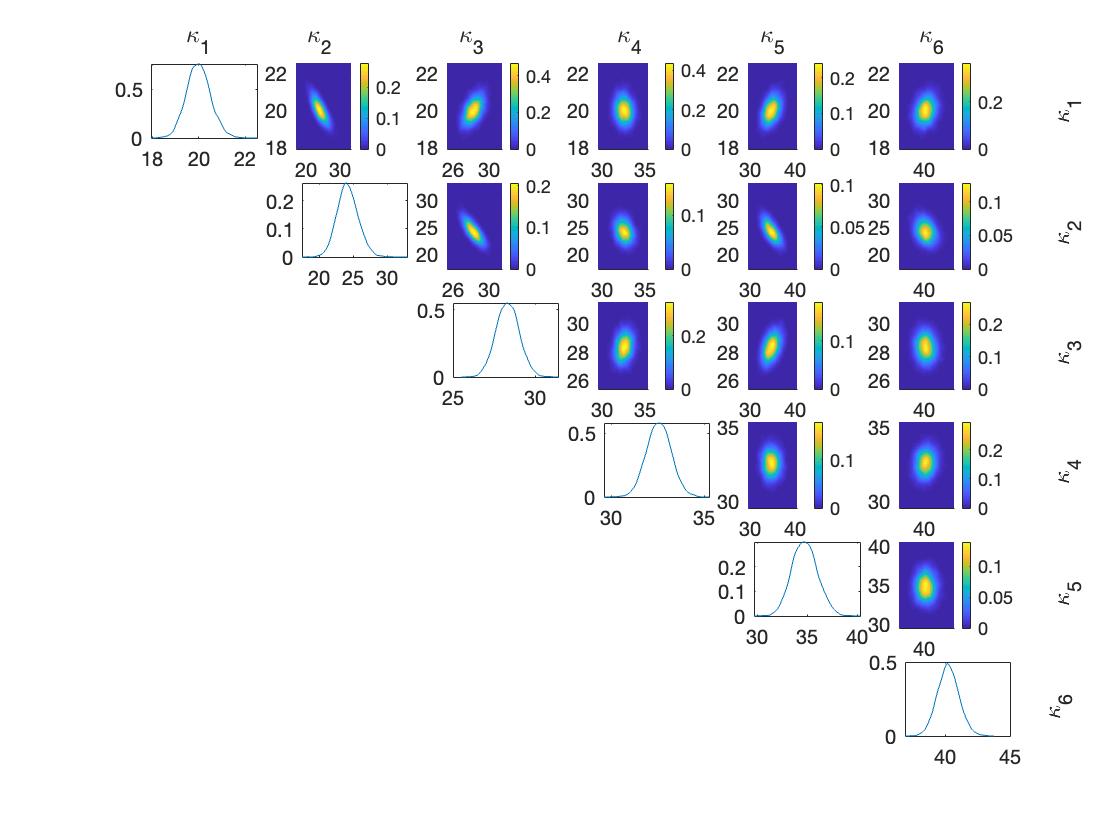}
\caption{The reconstructed posterior probability density for 6D IHCP.}
\label{fig4}

\end{figure}

 We also consider the IHCP with six unknown conductivities. For the corresponding setup one can refer to Fig. \ref{fig2}. The unknown conductivites $\boldsymbol{\kappa}=(\kappa_1, \cdots, \kappa_6)^\top$ are inferred with $N=20$ noisy measurements
  $\boldsymbol{\tilde{y}}=(u(\boldsymbol{x}_1), \cdots, u(\boldsymbol{x}_{20}))^\top$ at the measurements $(\boldsymbol{x}_1, \boldsymbol{x}_2, \cdots, \boldsymbol{x}_{20})^\top$. The noise $\delta$ is taken as $0.05$. We display the reconstructed probability in Fig. \ref{fig4}.


 \subsection{Inverse acoustic obstacle scattering}
 
 We apply the proposed algorithm to an inverse acoustic scattering problem with a sound-soft obstacle. We consider the scattering by long cylindrical obstacles with cross sections $\Omega\subset\mathbb{R}^2$ 
 that are starlike with respect to the origin. 
 In mathematics, we assume that $\Omega\subset\mathbb{R}^2$ is a bounded, 
 simply connected domain with $C^2$ boundary $\partial\Omega$. Then $\partial\Omega$ can be uniquely represented by a periodic function $r: [0, 2\pi)\rightarrow\mathbb{R}^+$:
 \begin{align}\label{parq1}
 \partial\Omega:=r(s)(\cos s, \sin s)=\exp(q(s))(\cos s, \sin s), \,\, s\in[0, 2\pi),
 \end{align}
 where $q(s)=\log r(s)$, $0<r(s)<r_{\max}$. 
 
  For a given incident plane wave
 \begin{align}
 u^{\rm i}(\boldsymbol{x}):=\exp(\mathfrak{i}k\boldsymbol{x}\cdot \boldsymbol{d}), \,\, \boldsymbol{x}\in\mathbb{R}^2,\,\, \boldsymbol{d}\in\mathbb{S}:=\{\hat{\boldsymbol{x}}\in\mathbb{R}^2\big| |\hat{\boldsymbol{x}}|=1\},
 \end{align}
 where $k>0$ is the wavenumber, $\mathfrak{i}=\sqrt{-1}$ and $\boldsymbol{d}:=(\cos\varsigma, \sin\varsigma)$ is the direction, the scattering problem is to find the scattered field $u^{\rm s}$, or the total field $u=u^{\rm i}+u^{\rm s}$, such that
 \begin{align}\label{sca1}
 \begin{aligned}
 &\Delta u+k^2u=0, &\text{in}\,\, \mathbb{R}^2\backslash\bar{\Omega},\\
 &u=0, &\text{on}\,\, \partial\Omega,\\
 &\lim_{r\rightarrow\infty}\sqrt{r}\left(\frac{\partial u^{\rm s}}{\partial r}-\mathfrak{i}ku^{\rm s}\right)=0, & \text{Sommerfeld radiation condition}.
 \end{aligned}
 \end{align}
It is well-known from the Sommerfeld radiation condition that  the scattered field admits the following asymptotic expansion
\begin{align}
u^{\rm s}(\boldsymbol{x}, \boldsymbol{d})=\frac{\exp(\mathfrak{i}\frac{\pi}{4})}{\sqrt{8k\pi}}\frac{\exp(\mathfrak{i}kr)}{\sqrt{r}}\left\{u^\infty(\boldsymbol{\hat{x}}, \boldsymbol{d})+O\left(\frac{1}{r}\right)\right\}\,\, \text{as}\,\, r:=|\boldsymbol{x}|\rightarrow\infty
\end{align}
 uniformly for all directions $\boldsymbol{\hat{x}}=\boldsymbol{x}/|\boldsymbol{x}|$. The function $u^\infty(\boldsymbol{\hat{x}}, \boldsymbol{d})$  defined on the unit circle $\mathbb{S}\subset\mathbb{R}^2$ is called the far field pattern of $u^{\rm s}$. 

The inverse scattering problem considered in this paper is to determine $\partial\Omega$ from the observation of the far field pattern $u^\infty(\boldsymbol{\hat{x}}, \boldsymbol{d})$.
We first list the basic knowledges to establish the forward map $\mathcal{G}$. Recall that the fundamental solution $\varPhi(\boldsymbol{x}, \boldsymbol{\tilde{x}})$ of the Helmholtz equation is given by
\begin{align*}
\varPhi(\boldsymbol{x}, \boldsymbol{\tilde{x}})=\frac{\mathfrak{i}}{4}H_0^1(|\boldsymbol{x}-\boldsymbol{\tilde{x}}|),
\end{align*}
 where $H_0^1$ is the Hankel function of the first kind of order zero. The single-layer potential operator $\mathcal{S}$ and the double-layer potential operator $\mathcal{K}$ are defined by 
 \begin{align*}
 \begin{aligned}
 (\mathcal{S}\varphi)(\boldsymbol{x})=2\int_{\partial\Omega}\varPhi(\boldsymbol{x}, \boldsymbol{\tilde{x}}) \varphi(\boldsymbol{\tilde{x}})ds(\boldsymbol{\tilde{x}}),\,\, \boldsymbol{x}\in\partial\Omega,\\
 (\mathcal{K}\varphi)(\boldsymbol{x})=2\int_{\partial\Omega} \frac{\partial\varPhi(\boldsymbol{x}, \boldsymbol{\tilde{x}})}{\partial \nu(\boldsymbol{\tilde{x}})}\varphi(\boldsymbol{\tilde{x}})ds(\boldsymbol{\tilde{x}}),\,\, \boldsymbol{x}\in\partial\Omega,
 \end{aligned}
 \end{align*}
 respectively. From \cite{Colton}, we can know that $\mathcal{S}$ and $\mathcal{K}$ are bounded from $C^{0, \alpha}(\partial\Omega)$ into $C^{1, \alpha}(\partial\Omega)$, $\alpha\in (0, 1)$. According to the single- and double-layer potentials, the scattered 
 field can be written as
 \begin{align}
 u^{\rm s}(\boldsymbol{x}; \Omega)=\int_{\partial\Omega} \left\{\frac{\partial\varPhi(\boldsymbol{x}, \boldsymbol{\tilde{x}})}{\partial\nu(\boldsymbol{\tilde{x}})}-\mathfrak{i}\tau\varPhi(\boldsymbol{x}, \boldsymbol{\tilde{x}})\right\}\varphi(\boldsymbol{\tilde{x}})ds(\boldsymbol{\tilde{x}}), \,\, \boldsymbol{x}\in\mathbb{R}^2 \backslash \bar{\Omega},
 \end{align}
 where $\tau$ is a real coupling parameter and $\varphi(\boldsymbol{\tilde{x}})$ is the unknown density function. Then the direct scattering problem is to find the density $\varphi$ such that 
 \begin{align}\label{sd1}
 (\mathcal{I}+\mathcal{K}-\mathfrak{i}\tau\mathcal{S})\varphi=-2u^{\rm i}\,\, \text{on}\,\, \partial\Omega. 
 \end{align}
There exists a unique solution $\varphi$ satisfying \eqref{sd1} and depending continuously on $u^{\rm i}$ \cite{Colton}. Furthermore, the far field pattern has the following form
\begin{align}\label{sd2}
u^\infty(\boldsymbol{\hat{x}}, \boldsymbol{d})=\frac{\exp(-\mathfrak{i}\frac{\pi}{4})}{\sqrt{8\pi k}}\int_{\partial\Omega} \left(k\nu(\boldsymbol{\tilde{x}})\cdot \boldsymbol{\hat{x}}+\tau\right)\exp(-\mathfrak{i}k\boldsymbol{\hat{x}}\cdot \boldsymbol{\tilde{x}})\varphi(\boldsymbol{\tilde{x}})ds(\boldsymbol{\tilde{x}}).
\end{align} 
By combining \eqref{sd1} with \eqref{sd2}, the direct scattering problem can be written as 
\begin{align}
u^\infty(\boldsymbol{\hat{x}}, \boldsymbol{d})=\mathcal{G}(\Omega), 
\end{align}
 where $\mathcal{G}$ is the shape-to-measurement operator. Using the parameterization \eqref{parq1} and taking the noise in measurements into account, the inverse model is given by
 \begin{align}
 \boldsymbol{\tilde{y}}=\mathcal{G}(q)(\boldsymbol{\hat{x}}, \boldsymbol{d})+\Xi,\,\, (\boldsymbol{\hat{x}}, \boldsymbol{d})\in \Gamma^{\rm o}\times\Gamma^{\rm i}\subset \mathbb{S}\times\mathbb{S},
 \end{align}
 where $\Gamma^{\rm o}$ is the aperture of the observation and   $\Gamma^{\rm i}$ is the aperture of the incident wave. 
 
 The prior $q$ is taken as the truncated Fourier series \cite{Li}
  \begin{align}\label{tr1}
 q_N(s)=\frac{\kappa_0}{\sqrt{2\pi}}+\sum_{n=1}^N \left(\frac{\kappa_n}{n^{\upsilon}}\frac{\cos ns}{\sqrt{\pi}}+\frac{\tilde{\kappa}_n}{n^{\upsilon}}\frac{\sin ns}{\sqrt{\pi}}\right),
 \end{align}
 where $\kappa_n$ and $\tilde{\kappa}_n$ are i.i.d. (independent and identically distributed) with $\kappa_n, \tilde{\kappa}_n\sim N(0, 1)$ and $\upsilon$ is a positive constant.
The ansatz data is used to the numerical reconstruction. 
 The synthetic Gaussian noise is added to the true forward model, i.e.,
 \begin{align}
 \boldsymbol{\tilde{y}}=\mathcal{G}(q)(\boldsymbol{\hat{x}}, \boldsymbol{d})+\delta\|\mathcal{G}(q)(\boldsymbol{\hat{x}}, \boldsymbol{d})\|(\xi_1 +\xi_2 \mathfrak{i}),
 \end{align}
 where $\xi_1, \xi_2$ are i.i.d. normal Gaussian.
 
\begin{figure}
\centering
\subfigure[Threelobes with incident angles $\varsigma=\frac{\pi}{2}$ and $\varsigma=\frac{3\pi}{2}$.]{\includegraphics[width=4.8cm]{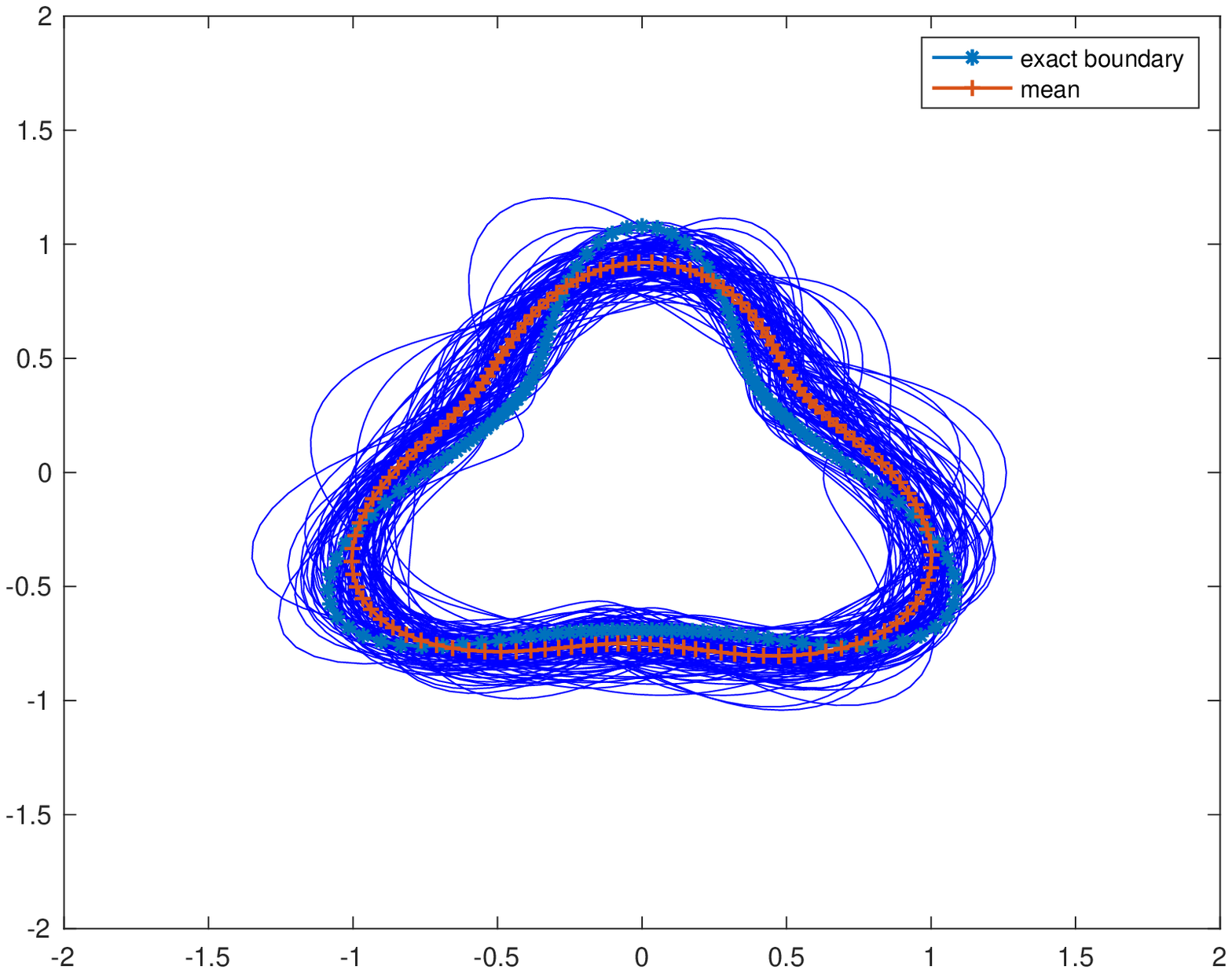}}
\subfigure[Pear with incident angles $\varsigma=0$ and $\varsigma=\frac{\pi}{2}$.]{\includegraphics[width=4.8cm]{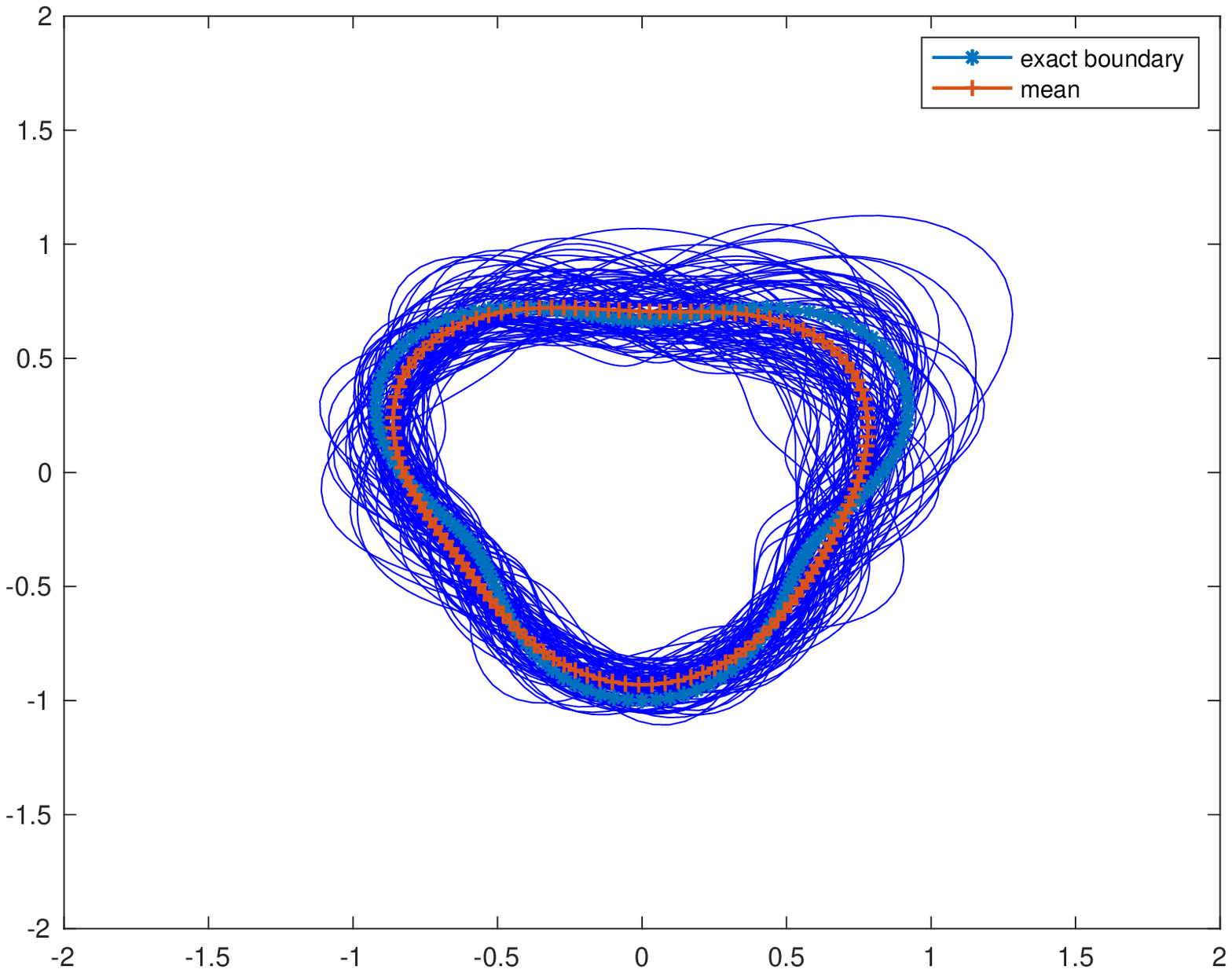}}\\
\subfigure[Bean with incident angles $\varsigma=0$ and $\varsigma=\frac{\pi}{2}$.]{\includegraphics[width=4.8cm]{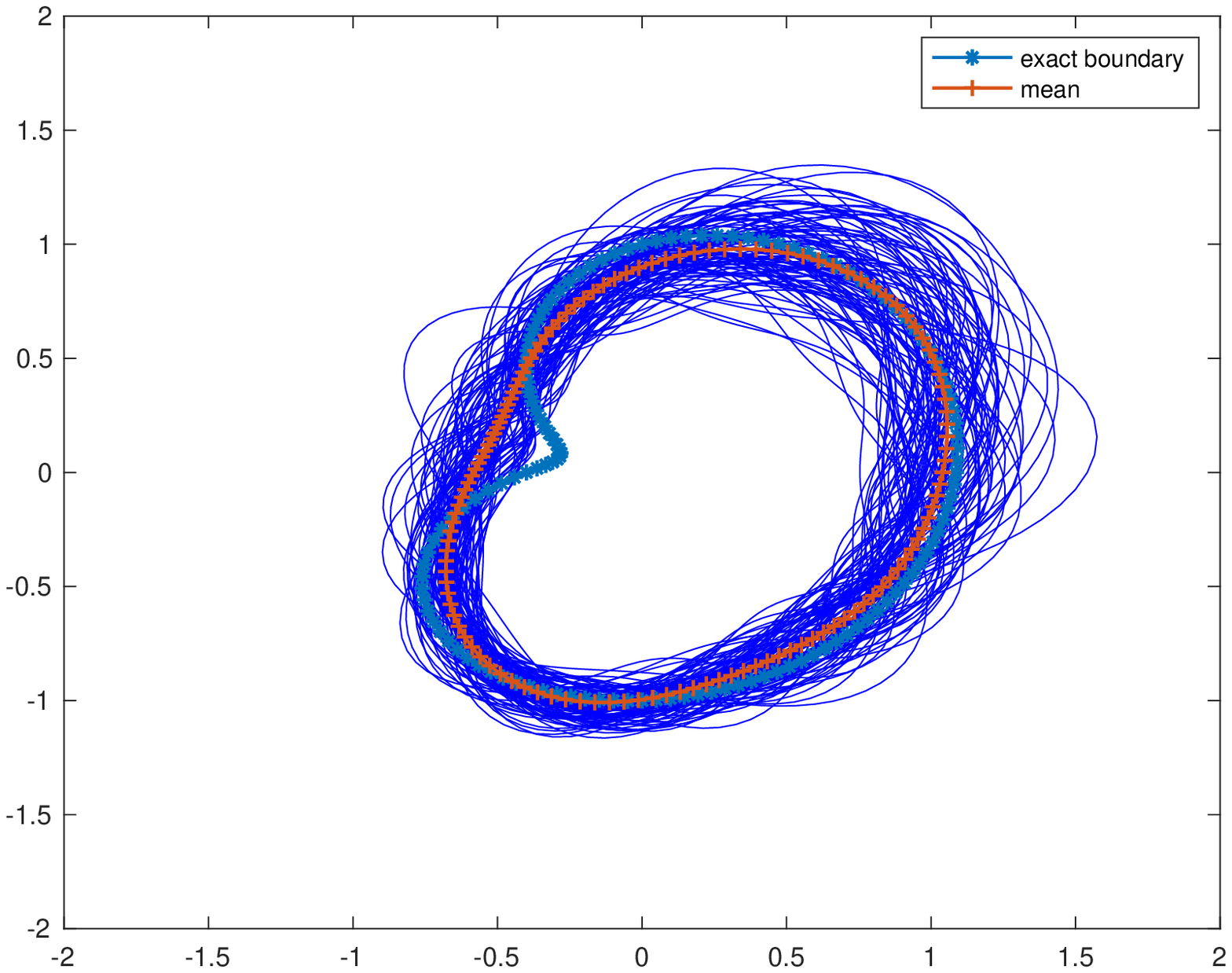}}
\subfigure[Peanut with incident angles $\varsigma=0$ and $\varsigma=\frac{\pi}{2}$.]{\includegraphics[width=4.8cm]{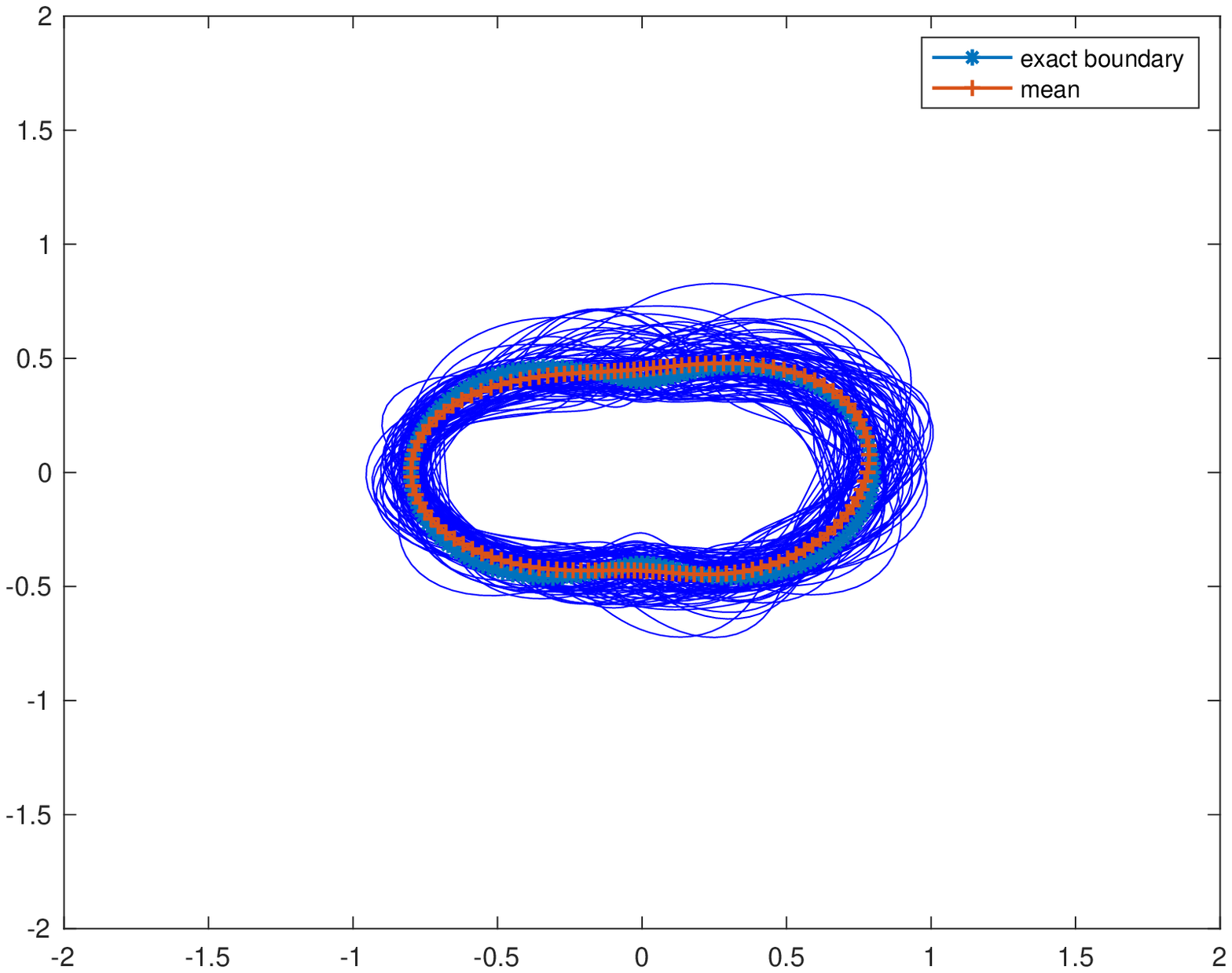}}\\
\subfigure[Acorn  with incident angles  $\varsigma=0$ and $\varsigma=\frac{\pi}{2}$.]{\includegraphics[width=4.8cm]{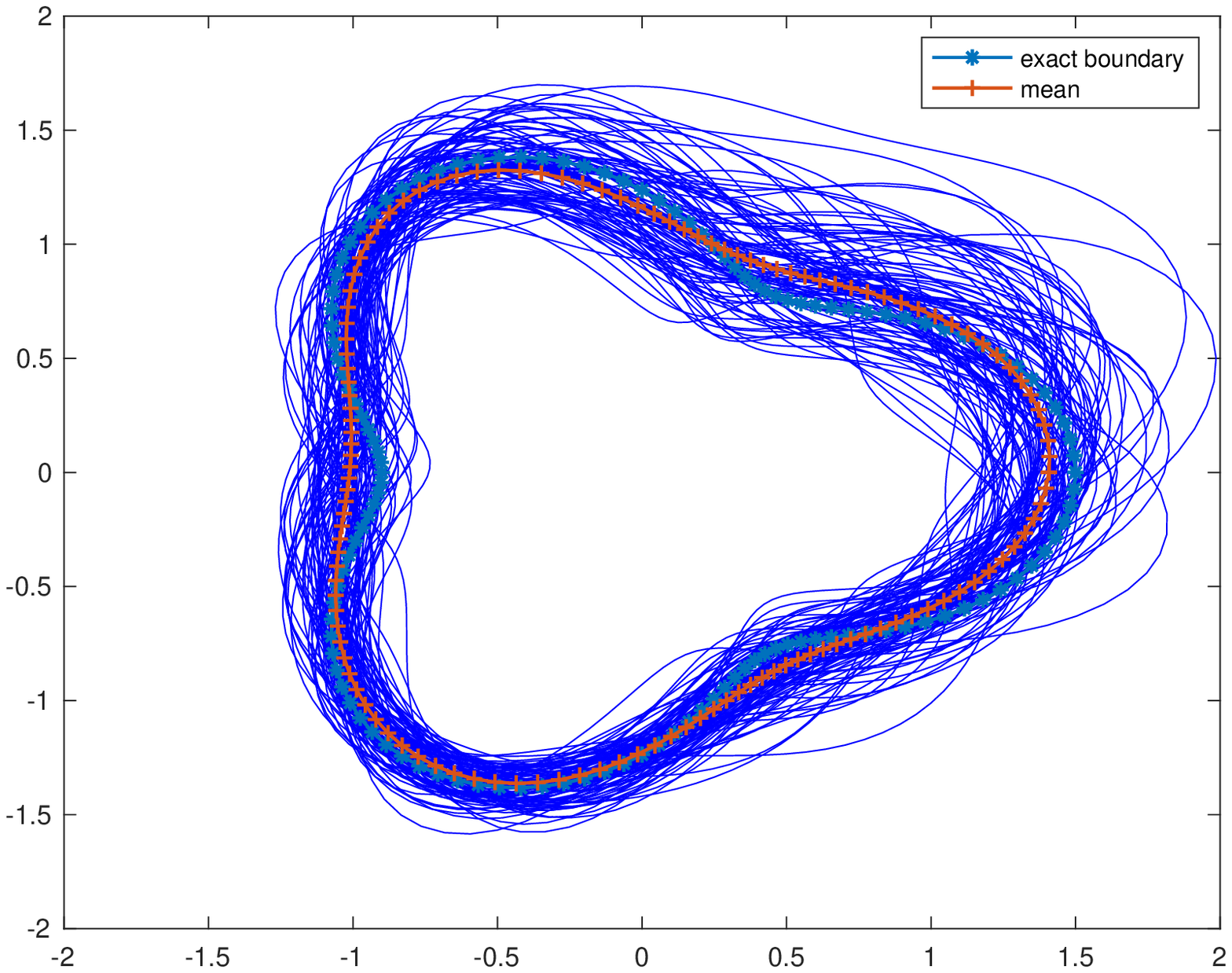}}
\subfigure[Roundedtriangle  with  incident angles $\varsigma=0$ and $\varsigma=\pi$.]{\includegraphics[width=4.8cm]{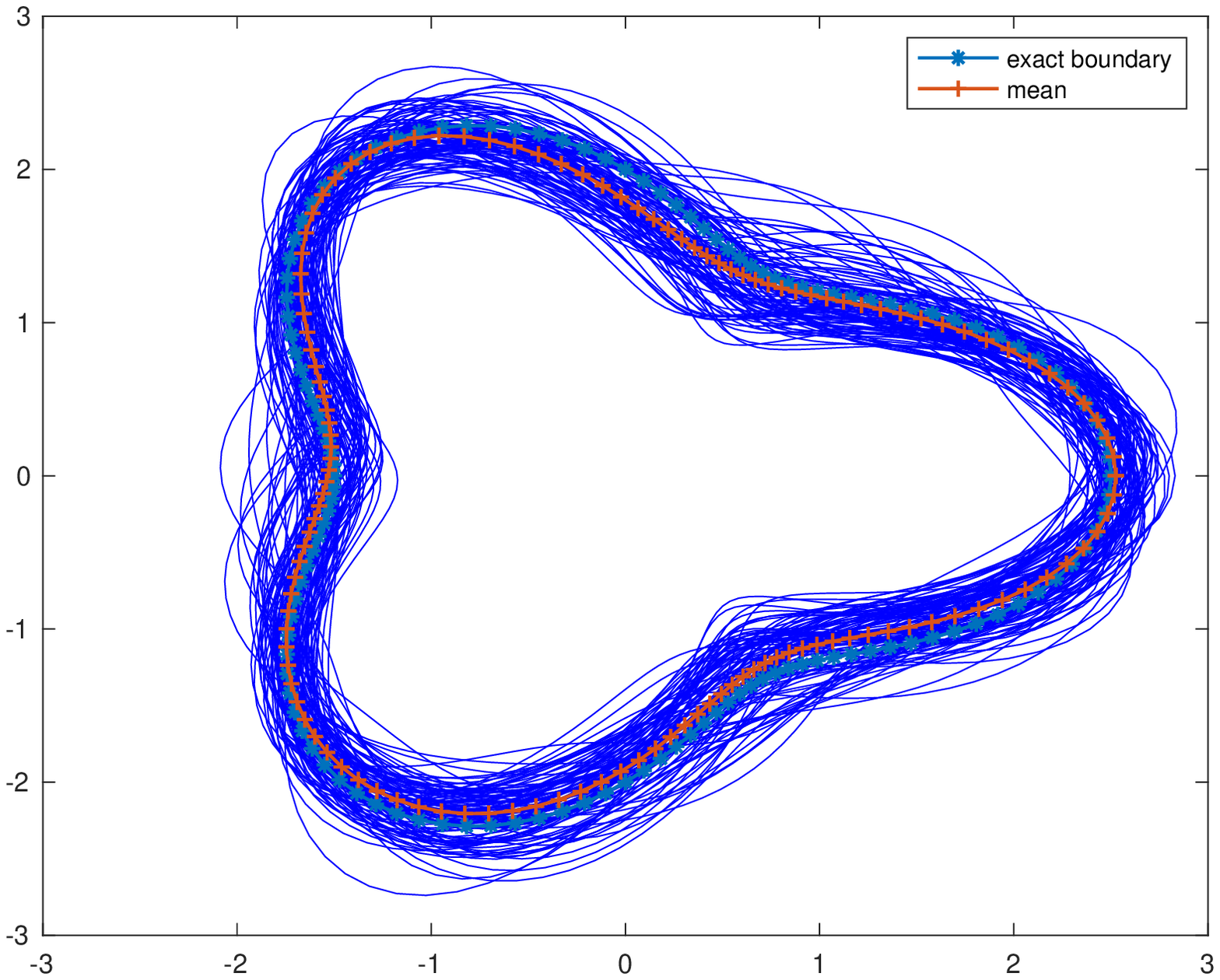}}\\
\subfigure[Roundrect  with incident angles $\varsigma=0$ and $\varsigma=\pi$.]{\includegraphics[width=4.8cm]{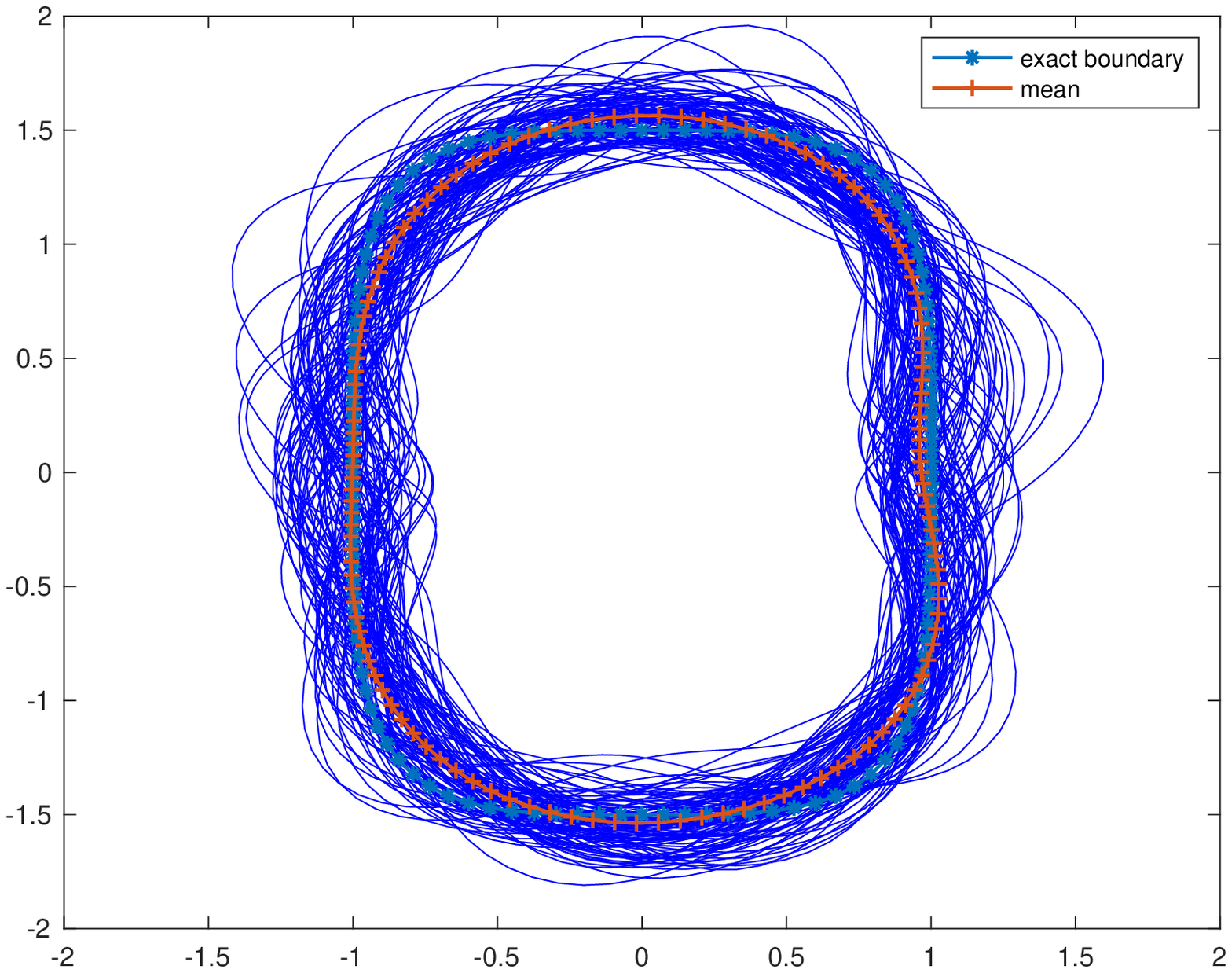}}
\subfigure[Kite  with incident angles $\varsigma=0$ and $\varsigma=\pi$.]{\includegraphics[width=4.8cm]{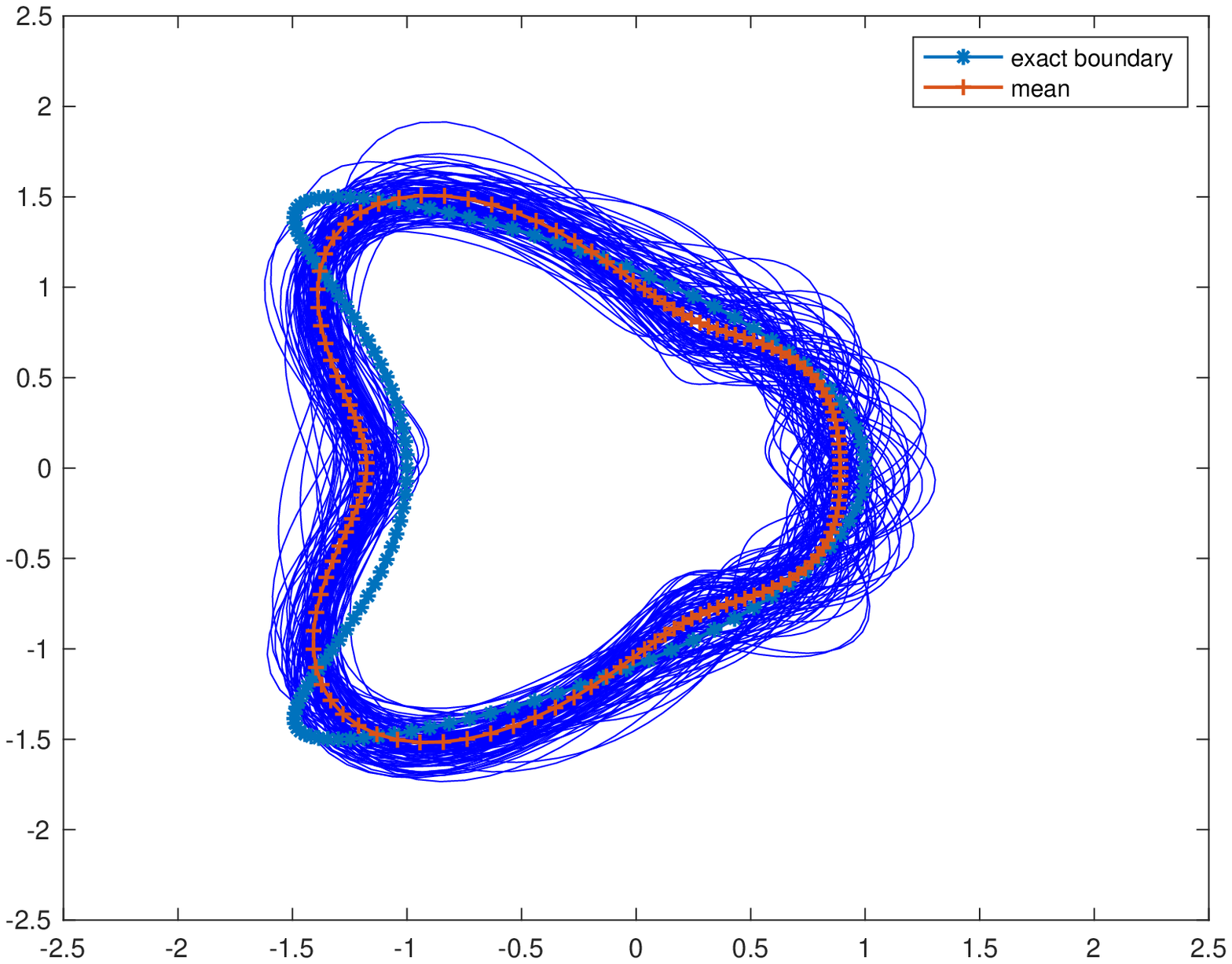}}
\caption{The numerical reconstructions for inverse scattering problem using full aperture data with two incident waves. The blue curves are the posterior samples drawn from the approximated posterior distribution.}
\label{fig5}

\end{figure}

 In the numerical implementation, we fix wavenumber $k=1$, $\delta=0.01$ and take $\upsilon=2.2, N=5$ in  \eqref{tr1}. Some frequently used test shapes in obstacle scattering are chosen as the test examples.
 For two incident waves, we collect the full aperture data. The reconstruction results are displayed in Fig. \ref{fig5}. 
\begin{itemize}[labelwidth={1em},font=\bfseries,align=left]
\item[(a)] threelobes: $r(s)=0.5+0.25\exp(-\sin 3 s)-0.1\sin s$;

\item[(b) ] pear: $r(s)=\frac{5+\sin 3s}{6}$;

\item[(c)] bean: $r(s)=\frac{1+0.9\cos s+0.1\sin 2s}{1+0.75\cos s}$;

\item[(d)] peanut: $r(s)=0.4\sqrt{4\cos^2 s+\sin^2 s}$;

\item[(e)] acorn: $r(s)=\frac{3}{5}\sqrt{\frac{17}{4}+2\cos 3s}$;
\item[(f) ] roundedtriangle: $r(s)=2+0.5\cos s$;

\item[(g) ] roundrect: $r(s)=( \cos^4 s+(2/3\sin s)^4 )^{-1/4}$;

\item[(h)] 
 kite: $x_1=\cos s+0.65\cos 2s-0.65$,     $ x_2= 1.5\sin s$.

\end{itemize}

 %

%

%

\bibliographystyle{plain}
\bibliography{ref}

\begin{thebibliography}{10}

\bibitem{Colton}
David Colton and Rainer Kress.
\newblock {\em Inverse acoustic and electromagnetic scattering theory}.
\newblock Springer, New York, 3rd edition, 2013.

\bibitem{Dick}
Josef Dick.
\newblock High-dimensional integration: The quasi-monte carlo way.
\newblock {\em Acta Numerica}, pages 133--288, 2013.

\bibitem{Gerstner}
Thomas Gerstner and Michael Griebel.
\newblock Numerical integration using sparse grids.
\newblock {\em Numerical Algorithms}, 18:209--232, 1998.

\bibitem{Hagmar}
Jonas Hagmar, Mats Jirstrand, Lennart Svensson, and Mark Morelande.
\newblock Optimal parameterization of posterior densities using homotopy.
\newblock In {\em 14th International Conference on Information Fusion}. IEEE,
  July 2011.

\bibitem{Hanebeck}
Uwe~D. Hanebeck, Kai Briechle, and Andreas Rauh.
\newblock Progressive bayes: A new framework for nonlinear state estimation.
\newblock {\em Proc. SPIE 2003, Multisource Information Fusion: Architectures,
  Algorithms, and Applications, B.V. Dasarathy}, 5099:256--267, 2003.

\bibitem{Li}
Zhaoxing Li, Zhiliang Deng, and Jiguang Sun.
\newblock Extended-sampling-bayesian method for limited aperture inverse
  scattering problems.
\newblock {\em SIAM J. IMAGING SCIENCES}, 13(1):422--444, 2020.

\bibitem{Nagel}
Joseph~B. Nagel and Bruno Sudret.
\newblock Spectral likelihood expansions for bayesian inference.
\newblock {\em Journal of Computational Physics}, (309):267--294, 2016.

\bibitem{Stuart}
Andrew~M. Stuart.
\newblock Inverse problems: A bayesian perspective.
\newblock {\em Acta Numerica}, 19:451--559, 2010.

\bibitem{Tarantola}
Albert Tarantola.
\newblock {\em Inverse Problem Theory and Methods for Model Parameter
  Estimation}.
\newblock Society for Industrial and Applied Mathematics, 2005.

\bibitem{Wagner}
Pau-Remo Wagner, Stefano Marelli, and Bruno Sudret.
\newblock Bayesian model inversion using stochastic spectral embedding.
\newblock {\em Journal of Computational Physics}, 436(110141), 2021.

\end{thebibliography}

\end{document}